\documentclass[a4paper,12pt]{amsart}
\usepackage[margin=29mm]{geometry}

\usepackage{amssymb}
\usepackage{amsmath}
\usepackage{amsthm}
\usepackage{amscd}
\usepackage[T1]{fontenc}
\usepackage{eucal,mathrsfs,dsfont}

\newcommand{\comp}{\mathds C}

\newcommand{\nat}{\mathds N}

\newcommand{\real}{\mathds R}
\newcommand{\rn}{{{\mathds R}^n}}
\newcommand{\rd}{{{\mathds R}^d}}
\newcommand{\Ee}{\mathds E}
\newcommand{\Pp}{\mathds P}
\newcommand{\I}{\mathds 1}
\newcommand{\Ff}{\mathcal{F}}

\renewcommand{\leq}{\leqslant}
\renewcommand{\geq}{\geqslant}
\renewcommand{\Re}{\ensuremath{\operatorname{Re}}}

\newtheorem{theorem}{Theorem}
\newtheorem{lem}[theorem]{Lemma}
\newtheorem{prop}[theorem]{Proposition}

\theoremstyle{definition}
\newtheorem{rem}[theorem]{Remark}

\newtheorem*{ack}{Acknowledgement}
\theoremstyle{remark}
\newtheorem{exa}{Example}

\begin{document}
\title[Transition probability densities of L\'evy processes]
{A note on the existence of transition probability densities for L\'evy processes}

\author{V.\ Knopova \and R.L.\ Schilling}\thanks{%
    \emph{R.L.\ Schilling}: Institut f\"ur Mathematische Stochastik, Technische
    Universit\"at Dresden, 01062 Dresden, Germany, \texttt{rene.schilling@tu-dresden.de}
}
\thanks{%
    \emph{V.\ Knopova}: V.M.Glushkov Institute of Cybernetics NAS of Ukraine,   03187, Kiev,
    Ukraine,
\texttt{vic$\,_{-}$knopova@gmx.de} }
\date{}

\maketitle

\begin{abstract}
    We prove several necessary and sufficient conditions for the existence of
    (smooth) transition probability densities for L\'evy processes and isotropic
    L\'evy processes. Under some mild conditions on the characteristic exponent we
    calculate the asymptotic behaviour of the transition density as $t\to 0$ and $t\to\infty$
    and show a ratio-limit theorem.



    \medskip\noindent\emph{Publication:} This paper has been published in the present version (but with different numbering) in \emph{Forum Mathematicum} \textbf{25} (2013) 125--149.
\end{abstract}

\section*{Introduction}
An $n$-dimensional L\'evy process is a stochastic process $X=(X_t)_{t\geq 0}$ with values in $\rn$, with independent and stationary increments, and with c\`adl\`ag (right-continuous, finite left limits) sample paths. It is well known that the transition probability $p_t$ of a L\'evy process can be characterized by the inverse Fourier transform:
\begin{equation*}
    \Ff^{-1}p_t(\xi) = \Ee_x e^{i X_t\cdot \xi} =e^{-t\psi(\xi)}, \quad t>0, \quad x, \xi \in\rn.
\end{equation*}
The function $\psi:\rn\to\comp$ is called the \emph{characteristic exponent} and it is determined by its \emph{L\'evy-Khintchine representation}
\begin{equation}\label{intro-e02}
    \psi(\xi)
    =i \ell\cdot\xi + \frac{1}{2} \,\xi\cdot Q\xi + \int_{\rn\setminus \{0\}} \left( 1-e^{iy\cdot\xi}+\frac{iy\cdot \xi}{1+|y|^2}\right)\nu(dy);
\end{equation}
here $\ell=(\ell^1,\ldots,\ell^n)\in\rn$, $Q = (q^{jk})\in\real^{n\times n}$ is a positive semi-definite matrix and $\nu$ is the L\'evy measure, i.e.\ a measure on $\rn\setminus\{0\}$ such that $\int_{\rn\setminus \{0\}} (1\wedge |y|^2)\,\nu(dy)<\infty$.  We call a L\'evy process $X$ \emph{symmetric} if for all $t\geq 0$ the random variables $X_t$ and $-X_t$ have the same distribution. Note that the characteristic exponent of symmetric L\'evy processes is real-valued and non-negative.

Many papers are devoted to distributional properties of L\'evy processes and to the existence of (necessary and) sufficient conditions under which the transition probability $p_t(dx)$ of a L\'evy process is absolutely continuous with respect to Lebesgue measure. The classic paper \cite{HW42} by Hartman and Wintner gives sufficient conditions in terms of the characteristic exponent $\psi$ under which there exists a transition density $p_t(x)$ of $X_t$; these conditions guarantee that $p_t \in C_\infty(\rn)$, where $C_\infty(\rn)$ denotes the set of all continuous functions which vanish at infinity. More precisely, if
\begin{equation}\label{hw}
    \lim_{|\xi|\to\infty} \frac{\Re \psi(\xi) }{\ln (1+|\xi|)} =\infty,
    \tag{$\textup{HW}_\infty$}
\end{equation}
then $p_t(dx) = p_t(x)\,dx$ for all $t>0$, and $p_t\in  L_1(\rn)\cap C_\infty(\rn)$. Also, if
\begin{equation}\label{hw1t}
    \liminf_{|\xi|\to\infty} \frac{\Re \psi(\xi) }{\ln (1+|\xi|)}>\frac{n}{t},
    \tag{$\textup{HW}_{1/t}$}
\end{equation}
then $p_s(dx)=p_s(x)\,dx$ for all $s\geq t$ and $p_s\in L_1(\rn)\cap C_\infty(\rn)$.

Note that the important issue is the speed at which the function $\psi$ tends to infinity. Since $\Ff^{-1} p_t = e^{-t\psi}$, the Riemann-Lebesgue lemma entails that
\begin{equation}\label{intro-e04}
    p_t(dx) = p_t(x)\,dx\text{\ \ for some $t>0$}
    \implies
    \lim_{|\xi|\to\infty} \Re\psi(\xi) = \infty,
\end{equation}
but this necessary condition does not tell anything about the rate of growth of $\psi$. Hartman and Wintner remark that \emph{the difficulties of the gap between these two conditions\textup{---i.e.\ \eqref{hw} and \eqref{intro-e04}---}are rather obscure} \cite[p.\ 287]{HW42}. Tucker \cite{Tu65} provides (complete but rather technical) necessary and sufficient criteria for the existence of a transition density. He mentions that Fourier analytic techniques \emph{are too crude for such a problem} (\cite[p.\ 317]{Tu65}). These two quotations are also the main motivation of our paper: to explain for which densities \eqref{hw} is indeed necessary and sufficient and how far we can get with Fourier analytic techniques.

Let us briefly review some of the other known criteria. Hawkes \cite{H79} shows that a L\'evy process has the strong Feller property---i.e.\ $x\mapsto \Ee_x f(X_t)$ is continuous for all bounded measurable functions $f$---if, and only if, the transition probabilities are absolutely continuous; in this case, the densities $p_t(x)$ are lower semicontinuous. Zabczyk \cite{Z70} shows that for an isotropic L\'evy process $X$ in $\rn$, $n\geq 2$, (see below for the precise definition) the following dichotomy holds: either $X$ is a compound Poisson process, or it is absolutely continuous for all $t>0$ with lower semi-continuous transition density.

Further sufficient conditions in dimension one were found by Kallenberg \cite[Section 5]{K81}:  if
\begin{equation}\label{kal}
    \lim_{\varepsilon \to 0} \frac{ \int_{-\varepsilon}^\varepsilon y^2 \nu(dy)}{\varepsilon^2\, |\ln \varepsilon|}=\infty,
    \tag{$\textup{K}_\infty$}
\end{equation}
then $p_t(dx) = p_t(x)\,dx$ for all $t>0$ and $p_t\in C_b^\infty (\real)\cap C_\infty(\real)$; and if
\begin{equation}\label{kal1t}
    \liminf_{\varepsilon \to 0} \frac{ \int_{-\varepsilon}^\varepsilon y^2 \nu(dy)}{\varepsilon^2\, |\ln \varepsilon|}> \frac 1t,
    \tag{$\textup{K}_{1/t}$}
\end{equation}
then  $p_s(dx)=p_s(x)\,dx$ for all $s\geq t>0$ and $p_s\in C_b^\infty (\real)\cap C_\infty(\real)$. For an $n$-dimensional analogue of \eqref{kal} and \eqref{kal1t} we refer to \cite{BK08}. See also  Orey \cite{O68} for yet another sufficient condition, as well as  the monograph of Sato \cite{S99} for more references and results.

\medskip
The main result of this note is to show for which class of L\'evy
processes the Hartman-Wintner condition \eqref{hw} is a necessary
and sufficient condition for the existence of a (smooth)
transition density $p_t(x)$. For isotropic processes we can
express \eqref{hw} in terms of the L\'evy measure $\nu$. Finally
we show that we can, under some mild conditions, express the
behaviour of the transition density at zero $p_t(0)$ in terms of
the measure of a ball with radius $t^{-1/2}$ in the metric given
by the characteristic exponent $\psi$. An important application of our results is a simple and at the same time sharp estimates for the on-diagonal behaviour or the transition probability density both as $t\to 0$ and $t\to\infty$, see Proposition \ref{app-21}. These estimates have an interpretation in geometric terms since it is given by the growth of balls related to a metric given by the characteristic exponent of the L\'evy process. These estimates extend previously known estimates for anisotropic stable and tempered stable processes, see \cite{St10}, \cite{R07},
\cite{St08}. In some sense, Theorem \ref{main-03} sharpens
\cite[Lemma 3.1]{St10} where an upper bound for the gradient of
the transition density $p_t$ of an $\alpha$-stable L\'evy process
is obtained  whose L\'evy measure $\nu$ is a $\gamma$-measure.
Although Theorem~\ref{main-03} does not provide a gradient estimate
in terms of the particular structure of the L\'evy measure, it
shows that the gradients of any order are in $L_1$.

\medskip\noindent
\textbf{Notation:} We denote by $\Ff u(\xi) = (2\pi)^{-n}\int_{\rn} u(x)\,e^{-ix\xi}\,dx$ the Fourier transform, $\Ff^{-1} w(x) = \int_{\rn} w(\xi)\,e^{i\xi x}\,d\xi$ is the inverse Fourier transform or characteristic function. By $J_\nu$ and $K_\nu$ we denote the (modified) Bessel functions of the first and third kind, cf.\ \cite{gra-ryz}. We write $C_\infty(\rn)$ for the continuous functions vanishing at infinity. Throughout this paper we will use the same letter $p_t$ to denote the transition probability $p_t(dx)$ and its density $p_t(x)$ w.r.t.\ Lebesgue measure. For functions $f(x)$ and $g(x)$ we write $f\asymp g$ if there are constants $c,C>0$ such that $cf(x)\leq g(x)\leq Cf(x)$ and we write $f\sim g$ (as $x\to a$) if $\lim_{x\to a} f(x)/g(x) = 1$. All other notation should be standard or self-explanatory.

\section{Main Results}
An $n$-dimensional (L\'evy) process $X$ is called \emph{isotropic}, if for any isometry $L:
\rn\to\rn$, $L(0)=0$, and all Borel sets $B\in\mathcal B(\rn)$
\begin{equation*}
    \Pp_0(X_t\in B)=\Pp_0(X_t\in LB), \quad t\geq 0.
\end{equation*}
In this case the L\'evy exponent is of the form $\psi(\xi)=g(|\xi|^2)$ for some continuous $g:[0,\infty)\to[0,\infty)$. For $n=1$ the notions of isotropy and symmetry coincide.

For an isotropic process we define $G(r):=-\omega_{n-1} \nu(B(0,r)^c)$, $n\geq 1$, where $\omega_{n-1}=2\pi^{n/2}/\Gamma\big(\frac n2\big)$ is a surface volume of the unit sphere $S^{n-1}\subset\rn$ and $\Gamma$ is Euler's Gamma function.


We can now state our main results.
\begin{theorem}\label{main-03}
    Let $X$ be an $n$-dimensional L\'evy process, $n\geq 1$, without a Gaussian component. The following conditions are equivalent:
    \begin{enumerate}
        \item[(a)] \eqref{hw};
        \item[(b)] for all $t>0$ the transition density exists, $p_t\in C^\infty(\rn)$ and  $\nabla^\alpha  p_t\in L_1(\rn)\cap C_\infty(\rn)$ for all  $\alpha\in\nat_0^n$;
        \item[(c)] for all $t>0$ the transition density exists and $p_t, \nabla p_t\in L_1(\rn)$.
    \end{enumerate}
    If $X$ is isotropic such that $\psi(\xi) = g(|\xi|^2)$ with an increasing function $g$, then the above conditions are equivalent to
    \begin{enumerate}
        \item[(d)] for all $t>0$ the transition density exists and $p_t\in C_\infty(\rn)$;
        \item[(e)] for all $t>0$ the transition density exists and $p_t\in L_\infty(\rn)$;
        \item[(f)] $e^{-t\psi}\in L_1(\rn)$ for all $t>0$.
    \end{enumerate}
\end{theorem}

\begin{theorem}\label{main-05}
    Let $X$ be an isotropic L\'evy process in $\rn$, $n\geq 2$, with characteristic exponent $\psi(\xi)=g(|\xi|^2)$.
    We assume that $X$ has no Gaussian component. Then \eqref{hw} is equivalent to
    \begin{equation}\label{main-e06}
        \lim_{\varepsilon\to 0} \frac{\nu(B(0,\varepsilon)^c)}{|\ln \varepsilon|}=\infty.
    \end{equation}
\end{theorem}

Before we proceed with the proofs of Theorem~\ref{main-03} and \ref{main-05} we add a few remarks and give some examples.

\begin{exa}\label{exa1}
    Since the existence and smoothness of a transition density is a time-dependent property, cf.\ \cite{S99}, the specification `\emph{for all $t>0$}' is essential in Theorem \ref{main-03}. A simple counterexample in dimension $n=1$ is the Gamma process, that is the L\'evy process with transition density
    $$
        p_t(x) = \frac{x^{t-1}}{\Gamma(t)}\,e^{-x},
        \quad t>0,\; x>0.
    $$
    It is not hard to see that its characteristic exponent is
    $$
        \psi(\xi) = \ln (1+i\xi)^{-1} = \frac 12\ln (1+\xi^2) + i\arctan(\xi).
    $$
    The two-sided (i.e.\ symmetrized) Gamma process whose transition density is $q_t := p_t*\tilde p_t$, $\tilde p_t(x) = p_t(-x)$, has
    \begin{align*}
        \psi(\xi)
        = \ln(1+\xi^2)
        &= \int_{\real\setminus\{0\}} \big(1-\cos(x\xi)\big) \left(\int_0^\infty \frac{1}{\sqrt{4\pi s}}\,e^{-s}\,e^{-|x|^2/(4s)}\,\frac{ds}s\right)\,dx\\
        &= \sqrt{\frac{2}{\pi}}\int_{\real\setminus\{0\}} \big(1-\cos(x\xi)\big) \frac{\exp(-|x|)}{\sqrt{|x|}}\,dx
    \end{align*}
    as characteristic exponent; this follows from a combination of \cite[9.23.4, 10.3]{ber-for} and \cite[8.437, p.\ 959]{gra-ryz}.

    Note that $p_t$ is for $t>1$ continuous and vanishes at infinity, while for $t=1$ the density is bounded and Borel measurable, while for $t\in (1/q,1)$ the density has a pole at $x=0$ but is still contained in $L_p$, $p$, $q$ being conjugate: $p^{-1} + q^{-1} = 1$. A similar picture is true for the density of the symmetrized process $q_t(x)$ which is given by $\Gamma(t)^{-1} \pi^{-1/2} (|x|/2)^{t-1/2} K_{t-1/2}(|x|)$ for $t>1/2$, cf.\
    \cite[17.344, p.\ 1151]{gra-ryz}.
\end{exa}

The following example explores the differences between the conditions \eqref{hw} and \eqref{kal}.
\begin{exa}\label{exa2}
 Let $n=1$. \eqref{kal} implies \eqref{hw}, but \eqref{hw} does not imply \eqref{kal}.
 Assume that $\ell=0, Q=0$ and $\nu(dy)= \frac 1{|y|}\, \ln\frac 1{|y|}\, \I_{B(0,1)}(y)dy$ in
 \eqref{intro-e02}. After some straightforward calculations we obtain that the related characteristic
 exponent behaves like $\psi(\xi)\sim \ln^2 |\xi |$ as $|\xi|\to\infty$  which implies \eqref{hw}. But
\begin{equation*}
    \lim_{\varepsilon\to 0} \frac{ 2\int_0^\varepsilon y\ln \frac{1}{y}\, dy}{\varepsilon^2\, |\ln \varepsilon|}
    =\lim_{\varepsilon\to 0} \frac{-2 \varepsilon \ln\varepsilon}{-2\varepsilon \ln \varepsilon -\varepsilon}
    =1,
\end{equation*}
    i.e.\  we have only \eqref{kal1t} for some $t>0$, but not \eqref{kal}.
\end{exa}

The next example shows that the condition $n\geq 2$ in \cite{Z70} is essential: for $n=1$ one can construct symmetric (hence, isotropic!) but continuous singular processes; the example shows also that the characteristic exponents can oscillate as $\xi\to\infty$.
\begin{exa}\label{exa3}
    Fix any $\lambda \in (0,2)$ and choose $M = M(\lambda) \in \mathds{N}$ such that $M > \frac{2}{2-\lambda}$. It is not hard to see that
    \begin{equation}\label{main-e08}
         \nu(dx) :=
         \sum_{j=1}^\infty \frac 12\,2^{\lambda M^j - j} \,\big(\delta_{2^{-M^j}}(dx) + \delta_{-2^{-M^j}}(dx)\big)
    \end{equation}
    is a L\'evy measure. The corresponding characteristic exponent is of the form
\begin{equation*}
     \psi(\xi)
     = \int_{\real\setminus\{0\}} (1-\cos (x \xi)) \nu(dx)
     = \sum_{j=1}^\infty 2^{\lambda M^j -j}\,\left( 1 - \cos \left(2^{-
     M^j}\xi\right)\right).
\end{equation*}
     This function enjoys the following properties:
\begin{align}
    \liminf_{|\xi|\to\infty} \psi(\xi)
    &=\label{main-e10}
    0;\\
    \limsup_{|\xi|\to\infty} \, \frac{\psi(\xi)}{|\xi|^{\lambda -\epsilon}}
    &=\label{main-e12}
    \infty \quad\text{for all}\quad\epsilon > 0;\\
    \lim_{|\xi| \to \infty} \, \frac{\psi(\xi)}{|\xi|^{\lambda + \epsilon}}
    &=\label{main-e14}
    0 \quad\text{for all}\quad\epsilon > 0.
\end{align}
The corresponding transition probability $p_t(dx)$ is not absolutely continuous, otherwise we would have by \cite{HW42} $\lim_{|\xi|\to \infty}\psi(\xi)=\infty$. Hence by \cite{HW42}, see also \cite[Theorem 27.16]{S99}, $p_t(dx)$ is continuous singular.

\medskip\noindent
\emph{Indeed:} We will need the following elementary inequalities
\begin{equation}\label{main-e16}
     \frac{t^2}2 \left( 1- \frac{t^2}{12}\right)^+
     \leq 1 - \cos(t)
     \leq \frac{t^2}2
\end{equation}
and
\begin{equation}\label{main-e18}
     1 - \cos(t)
     \leq 2^{1-\lambda}|t|^\lambda.
\end{equation}

Consider the sequence $2\pi 2^{M^k}$, $k\in\mathds N$. Clearly, $\lim_{k\to\infty} 2\pi 2^{M^k} = \infty$ and $2^{M^k} 2^{-M^j} \in \mathds N$ for $j\leq k$. Using \eqref{main-e16} and the fact that $\lambda < 2$ yields
\begin{align*}
     \psi\left(2\pi 2^{M^k}\right)
     &= \sum_{j=1}^\infty 2^{\lambda M^j - j} \left( 1 - \cos\left( 2\pi 2^{M^k}
     2^{-M^j}\right)\right)\\
     &= \sum_{j>k} 2^{\lambda M^j - j} \left( 1 - \cos\left( 2\pi 2^{M^k} 2^{-
     M^j}\right)\right)\\
     &\leq 2\pi^2 \sum_{j>k} 2^{(\lambda-2) M^j - j} \, 2^{2M^k}\\
     &\leq 2\pi^2 \sum_{j>k} 2^{(\lambda-2) M^{k+1} - j} \, 2^{2M^k}\\
     &= 2\pi^2 \, 2^{(M\lambda - 2M + 2)M^k}\sum_{j>k}2^{-j}.
\end{align*}
By the choice of $M = M(\lambda)$ we have $M\lambda -2M + 2 \leq 0$, thus $\lim_{k\to\infty} \psi\left( 2\pi 2^{M^k}\right) = 0$, and \eqref{main-e10} follows.

Using \eqref{main-e18} we find for any $\xi\in\rn$
\begin{align*}
    \psi(\xi)
    &= \sum_{j=1}^\infty 2^{\lambda M^j - j} \left( 1 - \cos\left( 2^{-M^j} \xi\right)\right)\\
    &\leq 2^{1-\lambda} \sum_{j=1}^\infty 2^{\lambda M^j - j} \, 2^{-\lambda M^j}\, |\xi|^\lambda\\
    &= 2^{1-\lambda}\,|\xi|^\lambda,
\end{align*}
which proves \eqref{main-e14}. For $\epsilon > 0$, $\xi\in\rn$ and any $k\in\mathds N$ we have by \eqref{main-e16}
\begin{align*}
    \frac{\psi(\xi)}{|\xi|^{\lambda - \epsilon}}
    &= |\xi|^{-\lambda + \epsilon} \sum_{j=1}^\infty 2^{\lambda M^j - j} \left( 1 - \cos \left( 2^{-M^j} \xi\right)\right)\\
    &\geq |\xi|^{-\lambda + \epsilon}  \, 2^{\lambda M^k - k} \left( 1 - \cos \left( 2^{-M^k} \xi\right)\right)\\
    &\geq \frac 12\,|\xi|^{2 -\lambda + \epsilon}  \, 2^{(\lambda-2) M^k - k} \left(
     1 - \frac{\xi^2}{12} \, 2^{-2M^k}\right)
\end{align*}
Setting $\xi_k := 2^{M^k}$ this gives
$$
    \frac{\psi(\xi_k)}{|\xi_k|^{\lambda - \epsilon}}
    \geq \frac12\, 2^{\epsilon M^k - k}\, \left(1 - \frac 1{12}\right)
    \xrightarrow{k\to\infty}\infty
$$
and \eqref{main-e12} follows.
\end{exa}

Lemma \ref{main-07} combined with Zabczyk's result \cite{Z70} indicates that \eqref{hw} is indeed not necessary for the existence of a $C_\infty$-transition density.
\begin{lem}\label{main-07}
    Let $X$ be an isotropic L\'evy process in $\rn$, $n\geq 1$.  Then
    \begin{equation}\label{main-e20}
    p_{t/2}(x) \text{\ \ exists and\ \ } p_{t/2} \in   L_1(\rn)\cap C_\infty(\rn)
    \implies
    \limsup_{|\xi|\to\infty} \frac{ \psi(\xi) }{\ln (1+|\xi|)}>\frac{n}{t}.
    \end{equation}
\end{lem}
\begin{proof}
Note that the left-hand side of \eqref{main-e20} implies $\|e^{-t\psi}\|_{L_1}<\infty$. Indeed, if $p_{t/2} \in C_\infty(\rn)\cap L_1(\rn)$, then $p_{t/2} \in L_\infty (\rn)\cap L_1(\rn)\subset L_2(\rn)$, hence $\Ff^{-1} p_{t/2} \in L_2(\rn)$. This shows that $|e^{-t\psi}|=|e^{-\frac t2 \psi}|^2 \in L_1(\rn)$.

Since $X$ is isotropic, $\psi(\xi)=g(|\xi|^2)$ for some continuous function $g:[0,\infty)\to[0,\infty)$, and we have
\begin{equation}\label{main-e22}
\begin{split}
    \|e^{-t\psi}\|_{L_1}
    &=\int_0^\infty m\left\{ \xi\in\rn:\,\, \psi(\xi)\leq -\frac{\ln s}{t}\right\} ds\\
    &=t\int_0^\infty m\left\{\xi\in\rn:\,\, g(|\xi|^2)\leq x\right\} e^{-tx}\, dx\\
    &\geq tc_n\int_0^\infty (g^{-1}(x))^{n/2}\, e^{-tx}\, dx,
\end{split}
\end{equation}
where $g^{-1}(x):=\inf\{u :\,\,g(u)\geq x\}$ and $m$ denotes
Lebesgue measure in $\rn$.  The left-hand side of \eqref{main-e22} is
finite and we conclude that $\liminf_{x\to\infty}
g^{-1}(x)e^{-2tx/n} = 0$.

By the very definition of the generalized inverse we get  $g^{-1}(g(x))\leq x \leq g(g^{-1})(x)$. Therefore,
\begin{equation}\label{main-e24}
    g^{-1}(x) e^{-2tx/n} \geq g^{-1}(x) e^{-2t g(g^{-1}(x))/n}\geq 0.
\end{equation}
Since the transition density exists, it is an easy consequence of the Riemann-Lebesgue lemma that $\lim_{|\xi|\to\infty}\psi(\xi)=\infty$, cf.\ \eqref{intro-e04}. Therefore, $g$ and $g^{-1}$ are onto. If we combine this fact with the inequality \eqref{main-e24}, we see that $\liminf_{x\to\infty} g^{-1}(x)e^{-2tx/n} = 0$ entails $\liminf_{u\to\infty} u e^{-2tg(u)/n}= 0$. Consequently,
\begin{equation*}
    \liminf_{u\to\infty}  \left( 1-\frac{2tg(u)}{n\ln u}\right) \ln u <0,
\end{equation*}
implying $\limsup_{u\to \infty} g(u^2) /\ln u>n/t$.
\end{proof}

A small modification of Example \ref{exa3} shows that we can have a density for all $t>0$ even if \eqref{hw} fails.
\begin{exa}\label{exa4}
    Let $\psi(\xi)$ be the characteristic exponent constructed in Example \ref{exa3} and observe that $\psi_1(\xi) := \psi(\xi) + \ln(1+|\xi|^2)$ is again a characteristic exponent of a L\'evy process $Y$. Since the L\'evy process corresponding to $\log(1+|\xi|^2)$ has a density for all $t>0$, cf.\ Example \ref{exa1}, we infer that all $Y_t$ are absolutely continuous.

    On the other hand,
    $$
        \limsup_{|\xi|\to\infty} \frac{\psi_1(\xi)}{\ln|\xi|} = \infty
        \quad\text{while}\quad
        \liminf_{|\xi|\to\infty} \frac{\psi_1(\xi)}{\ln|\xi|} = 2.
    $$
\end{exa}

%
%
\section{Proofs}
Let us now turn to the proof of  Theorem~\ref{main-03}. In a first step, Lemma \ref{proofs-09} below, we show that \eqref{hw} is also a sufficient condition for the existence of a $C^\infty_b \cap
C_\infty$ density.
\begin{lem}\label{proofs-09}
    Suppose that \eqref{hw} holds true. Then $p_t\in C^\infty(\rn)$, and  $\nabla^\alpha p_t\in L_2(\rn)\cap C_\infty(\rn)$ for all  $\alpha\in\nat_0^n$.
\end{lem}
\begin{proof}
Observe that \eqref{hw} implies $\frac{\Re\,\psi(\xi)}{\ln |\xi|}>c$ for any $c>0$ and all $\xi$ with  $|\xi|\geq |\xi_0(c)|$. For every $t>0$ and $k\geq 1$ we can achieve that $tc-k > n$ and that for $|\xi|$ large enough we get
\begin{align*}
    |\xi|^k \exp\big(-t\Re\, \psi(\xi)\big)
    =\exp\left(-\ln |\xi| \left[t\,\frac{\Re\, \psi(\xi)}{\ln |\xi|} - k\right]\right)
    \leq \exp\left(-(tc-k)\ln |\xi|\right).
\end{align*}
This shows
$$
    |\xi|^k e^{-t\psi(\xi)}\in L_2(\rn) \quad\text{for all\ } k\geq 1,
$$
which means that
\begin{equation*}
     \nabla^\alpha  p_t \in \bigcap_{k\geq 1} H^k(\rn)\hookrightarrow C_\infty(\rn)\quad\text{for all\ } \alpha\in\nat_0^n,
\end{equation*}
where $H^k(\rn)$ is the $L_2$-Sobolev space of order $k$.
\end{proof}

\begin{proof}[Proof of Theorem~\ref{main-03}]
(a)$\Rightarrow $(b): We decompose the characteristic exponent into a sum, $\psi = \psi_1 + \widetilde\psi_1$, where
\begin{align*}
    \psi_1(\xi)
    &= \int_{0<|y|\leq 1} \Big(1-e^{i y\cdot\xi} + iy\cdot\xi\Big)\,\nu(dy)\\
    \widetilde\psi_1(\xi)
    &= \int_{|y|> 1} \left(1-e^{iy\cdot\xi} + \frac{iy\cdot\xi}{1+|y|^2}\right)\nu(dy) + i\left(\ell\cdot\xi - \int_{0<|y|\leq 1} \frac{y\cdot\xi\,|y|^2}{1+|y|^2}\,\nu(dy)\right).
\end{align*}
By construction, $\psi_1$ and $\widetilde\psi_1$ are characteristic exponents of two independent L\'evy processes. Denote their transition probabilities by $p_{1,t}$ and $\widetilde p_{1,t}$, respectively. Because of independence, $p_t = p_{1,t}*\widetilde p_{1,t}$. Moreover, $\psi_1$ is infinitely often differentiable; for any multiindex $\alpha\in\nat^n$ we have
\begin{equation}\label{proofs-e26}
    \nabla^\alpha \psi_1(\xi)
    =
    \begin{cases}
    \displaystyle
        i\int_{0<|y|\leq 1} \big(1-e^{i y\cdot\xi}\big)\, y_\ell\,\nu(dy), & \alpha = e_\ell; \\[\bigskipamount]
    \displaystyle
         i^{|\alpha|-2}\int_{0<|y|\leq 1} e^{i y\cdot\xi}\,  y^\alpha\,\nu(dy), & |\alpha| \geq 2.
    \end{cases}
\end{equation}
(As usual, $y^\alpha := \prod_{j=1}^n y_j^{\alpha_j}$ for $\xi\in\rn$ and $\alpha\in\nat_0^n$)
This shows, in particular, that all derivatives of $\psi_1$ are polynomially bounded. Since $\Re\widetilde\psi_1$ is bounded,
$$
    \Re\psi_1(\xi) \leq \Re\psi(\xi) \leq \Re\psi_1(\xi) + \|\Re\widetilde\psi_1\|_\infty,
$$
indicating that \eqref{hw} holds for $\psi$ if, and only if, \eqref{hw} holds for $\psi_1$.

By the chain rule, for each $\alpha\in\nat^n$ there is a polynomial $P$ such that
$$
    \nabla^\alpha e^{-t\psi_1(\xi)}
    = (-t)^{|\alpha|}P\big(\psi_1,\ldots,\nabla^\alpha\psi_1\big) e^{-t\psi_1};
$$
since $\psi_1$ and all of its derivatives are polynomially bounded, \eqref{hw} implies that $\sup_\xi\big|\xi^\beta \nabla^\alpha e^{-t\psi_1(\xi)}\big|\leq c_{\alpha,\beta,t}<\infty$ for all $\alpha,\beta\in\nat_0^n$. Therefore, $\exp(-t\psi_1)\in S(\rn)$ where $S(\rn)$ denotes the Schwartz space of rapidly decreasing functions. Thus,
$$
    p_{1,t}(x) = \Ff e^{-t\psi_1}(x),
$$
i.e.\ $p_{1,t}(dx)=p_{1,t}(x)\,dx$ with a density from $S(\rn)$. The identity $p_t = p_{1,t}*\widetilde p_{1,t}$ shows that $p_t(dx)=p_t(x)\,dx$ with a $C^\infty$-density which satisfies
$$
    \nabla^\alpha p_t(x) = (\nabla^\alpha p_{1,t})*\widetilde p_{1,t}(x) = \int_\rn \nabla^\alpha  p_{1,t}(x-y)\,\widetilde p_{1,t}(dy).
$$
Using Fubini's theorem, and the fact that $\widetilde p_{1,t}$ is a probability measure, yields
$$
    \|\nabla^\alpha p_t\|_{L_1}
    \leq \iint |\nabla^\alpha p_{1,t}(x-y)|\,dx\,\widetilde p_{1,t}(dy)
    = \|\nabla^\alpha p_{1,t}\|_{L_1}.
$$
That $p_t$ and $\nabla^\alpha p_t$ are in $C_\infty$ follows from Lemma \ref{proofs-09}.

\medskip\noindent
(b)$\Rightarrow$(c): This is obvious.

\medskip\noindent
(c)$\Rightarrow$(a): Since $\nabla p_t\in L_1(\rn)$, the Riemann-Lebesgue lemma shows that
$$
    |\xi| e^{-t\psi(\xi)}
    = \exp\left(-\ln |\xi| \left[\frac{\psi(\xi)}{\frac 1t\,\ln |\xi|} - 1\right]\right)\in C_\infty(\rn)
$$
for all $t>0$. Letting $t\to 0$ implies \eqref{hw}.

\medskip
From now on we assume that $X$ is isotropic with $\psi(\xi)=g(|\xi|^2)$ and with an increasing function $g:[0,\infty)\to[0,\infty)$.

\medskip\noindent
(a)$\Rightarrow$(d)$\Rightarrow$(e): This follows from the above statements.

\medskip\noindent
(e)$\Rightarrow$(f): By assumption, $p_t\in L_\infty(\rn)\cap L_1(\rn)$ for all $t>0$. In particular $p_t\in L_2(\rn)$ and, by Plancherel's theorem, $\Ff^{-1} p_t = e^{-t\psi} \in L_2(\rn)$. Since this holds for all $t>0$, we see for $t=2s$ that $e^{-s\psi}\in L_1(\rn)$ for all $s>0$.

\medskip\noindent
(f)$\Rightarrow$(a): Since $e^{-t\psi}\in L_1$, we get $(2\pi)^n\,p_t(0) = \int_{\rn} e^{-t\psi(\xi)}\,d\xi<\infty$. Introducing polar coordinates and integrating by parts yields
\begin{align*}
    \int_{\rn} e^{-t\psi(\xi)}\,d\xi
    &= \int_{\rn} e^{-t g(|\xi|^2)}\,d\xi\\
    &= \omega_{n-1} \int_0^\infty e^{-t g(r^2)}\,r^{n-1}\,dr\\
    &\geq \omega_{n-1} \int_1^\infty e^{-t g(r^2)}\,r^{n-1}\,dr\\
    &= \frac{\omega_{n-1}}n\left( \lim_{s\to\infty} e^{-t g(s^2)}\,s^n - e^{-t g(1)}-
      \int_1^\infty r^n\,d_re^{-t g(r^2)}\right).
\end{align*}
Since $r\mapsto e^{-tg(r^2)}$ is decreasing, the integral appearing in the last line is negative and the calculation shows that $\lim_{s\to\infty} e^{-t g(s^2)}\,s^n$ is finite. Therefore, $e^{-tg(r^2)} \leq c_t\,r^{-n}$ for all $r>1$ and with some suitable constant $c_t<\infty$. Then
$$
    \frac{\psi(\xi)}{\ln |\xi|}
    = \frac{g(|\xi|^2)}{\ln |\xi|}
    \geq \frac{n}{t}-\frac{\ln c_t}{t \ln|\xi|},\quad |\xi| > 1,
$$
implying
$$
    \liminf_{|\xi|\to \infty} \frac{\psi(\xi)}{\ln|\xi|}\geq \frac{n}{t}.
$$
Letting $t\to 0$, we get $(HW_\infty)$.
\end{proof}

Let us now turn to the proof of Theorem \ref{main-05}. Recall that the Bessel function of the first kind, $J_\nu(z)$, is defined by
\begin{equation}\label{proofs-e28}
    J_\nu(z)
    =\sum_{n=0}^\infty \frac{(-1)^n}{\Gamma(n+\nu+1)n!} \left(\frac{z}{2}\right)^{2n+\nu}, \quad \nu,\, z\in\real.
\end{equation}

\begin{proof}[Proof of Theorem \ref{main-05}]
    Recall that $G(r) = -\omega_{n-1}\nu(B(0,r)^c)$, $\omega_{n-1} = 2\pi^{n/2}/\Gamma\big(\frac n2\big)$. Switching to polar coordinates in the L\'evy-Khintchine formula \eqref{intro-e02}, we get
    \begin{equation}\label{proofs-e30}
        g(u^2)=\int_0^\infty \big(1-H_{\frac{n-2}{2}} (ur)\big) \, dG(r),\quad n\geq 1,
    \end{equation}
    cf.\ \cite[p.\ 99]{B55}, where  $H_\nu(r) := 2^\nu\,\Gamma(\nu+1) r^{-\nu}J_\nu(r)$. Note that $H_\nu(0)=1$.  For $\varepsilon>0$ and $\gamma>1$ (which will be determined later) we split the integral expression appearing in \eqref{proofs-e30} into two parts and get
 \begin{align*}
    \frac{g(\varepsilon^{-2})}{|\ln \varepsilon|}
    =\frac{1}{|\ln \varepsilon|}\left( \int_0^{\gamma\varepsilon} +\int_{\gamma\varepsilon}^\infty \right) \big(1-H_{\frac{n-2}{2}}(r \varepsilon^{-1})\big)\,dG(r)
    =: I_1(\varepsilon)+I_2(\varepsilon).
\end{align*}
Note that both $I_1$ and $I_2$ are nonnegative. By \eqref{proofs-e28},
    $$
    r^{-\nu} J_\nu(r)=\frac{1}{2^\nu \Gamma(\nu+1)} \left[1-\frac{r^2}{4(\nu+1)}+o(r^2)\right]
    \quad \text{ as } \quad r\to\infty,
    $$
      implying that  there exist $c_1, c_2>0$ such that
\begin{equation*}
    c_1 r^2 \leq 1-H_{\frac{n-2}{2}}(r)\leq c_2 r^2
    \quad\text{for all\ \ } 0\leq r\leq \gamma,
\end{equation*}
which gives
\begin{equation*}
    c_1 \int_0^{\varepsilon\gamma} \frac{r^2}{\varepsilon^2\, |\ln \varepsilon |} \,dG(r)
    \leq I_1(\varepsilon)
    \leq c_2 \int_0^{\varepsilon \gamma}\frac{r^2}{\varepsilon^2\, |\ln \varepsilon |}\,dG(r).
\end{equation*}
Therefore
\begin{gather*}
    \lim_{\varepsilon\to 0} I_1(\varepsilon)=\infty
    \quad\text{if, and only if,}\quad
    \lim_{\varepsilon\to 0} \int_0^\varepsilon \frac{r^2}{\varepsilon^2\, |\ln \varepsilon |}\,dG(r)=\infty.
\end{gather*}

\noindent Consider the second term. From \cite[pp.\ 359, 368]{WW58} we know that
$$
    H_\nu(z)\sim \sqrt{\frac{2}{\pi}}\frac{1}{z^{\nu+\frac 12}}  \cos\Big(z-\frac{\pi \nu}{2}-\frac{\pi}{4}\Big), \quad z\to \infty.
$$
Therefore, there exist constants $C_1, C_2>0$ such that
\begin{equation}\label{proofs-e32}
    1-\frac{C_1}{r^{\frac{n-1}{2}}}
    \leq 1-H_{\frac{n-2}{2}}(r)
    \leq 1+ \frac{C_2}{r^{\frac{n-1}{2}}}
    \quad\text{for all\ \ } r\geq \gamma.
\end{equation}
If we replace in \eqref{proofs-e32} $r$ by $r\varepsilon^{-1}$ we get
\begin{equation*}
    1-C_1 \left(\frac{\varepsilon}{r}\right)^{\frac{n-1}{2}}
    \leq 1-H_{\frac{n-2}{2}}(r \varepsilon^{-1})
    \leq 1+ C_2 \left(\frac{\varepsilon}{r}\right)^{\frac{n-1}{2}};
\end{equation*}
for $r>\gamma\varepsilon$ with a sufficiently large constant $\gamma$ (depending only on $C_1, C_2$ and $n$), we get new constants $0<C_3, C_4<\infty$ such that
$$
    0 < C_3 \leq 1-H_{\frac{n-2}{2}}(r\varepsilon^{-1}) \leq C_4\quad\text{for all\ \ }r\geq \gamma\varepsilon.
$$
Integrating this expression over $[\gamma\varepsilon,\infty)$ w.r.t.\ $dG(r)$ reveals that $\frac{|G(\gamma\varepsilon)|}{|\ln \varepsilon|} \asymp I_2(\varepsilon)$. This shows
\begin{gather}
    \lim_{\varepsilon\to 0} I_2(\varepsilon)=\infty
    \quad\text{if, and only if,}\quad
    \lim_{\varepsilon\to 0} \frac{|G(\varepsilon)|}{|\ln \varepsilon|}=\infty.\label{proofs-e34}
\end{gather}
From  \eqref{proofs-e34} we see that   \eqref{main-e06} implies \eqref{hw}.

Conversely, assume that   \eqref{main-e06} fails, but that \eqref{hw} holds.  In this case there exists a sequence $\epsilon_k>0$ such that $\lim_{k\to\infty}\epsilon_k = 0$ and
$$
    \lim_{k\to \infty}  \left|\frac{G(\epsilon_k)}{\ln\epsilon_k}\right| < \infty.
$$
Integration by parts yields
$$
    \int_0^\epsilon r^2\,dG(r) = \epsilon^2G(\epsilon) - 2\int_0^\epsilon rG(r)\,dr.
$$
By l'Hospital's rule we see
\begin{align*}
    \lim_{k\to\infty} \frac{2\int_0^{\epsilon_k} rG(r)\,dr}{\epsilon_k^2 \ln\epsilon_k}
    = \lim_{k\to\infty} \frac{2\epsilon_k G(\epsilon_k)}{2\epsilon_k \ln\epsilon_k + \epsilon_k}
    = \lim_{k\to\infty} \frac{2G(\epsilon_k)}{2\ln\epsilon_k + 1}.
\end{align*}
Since the latter is finite, this shows that
$$
    \lim_{k\to\infty} \frac{\epsilon_k^2 G(\epsilon_k) - 2\int_0^{\epsilon_k} r G(r)\,dr}{\epsilon_k^2 |\ln\epsilon_k|} = 0
$$
and, in particular, $\lim_{k\to\infty} I_1(\epsilon_k) =0$. By construction, $\lim_{k\to\infty} I_2(\epsilon_k)<\infty$ and, therefore,
$$
    \lim_{k\to\infty} \frac{g(\epsilon_k^{-2})}{|\ln\epsilon_k|} < \infty,
$$
i.e.\ \eqref{hw} fails and we have reached a contradiction. The proof is complete.
\end{proof}

\section{Extensions}
Let us give two generalizations of Theorem \ref{main-03}. For this we need to recall the notions of de- and increasing rearrangements of a function. Our standard reference is the monograph \cite{ben-sha}. Let $u$ be a real-valued measurable function defined on a measurable subset $B\subset\rn$. As usual, $\mu_u(t):=m\{\xi\in B\::\: |u(\xi)|>t\}$ ($m$ is Lebesgue measure) is the \emph{distribution function} of $u$ and
\begin{equation*}
    u^*(s):=\inf\{t\geq 0\::\: \mu_u(t)\leq s\}=\sup\{t\geq 0\::\: \mu_u(t) > s\}
\end{equation*}
is called the \emph{decreasing rearrangement} of $u$. Note that $u^* : [0,\infty)\to [0,\infty]$ is decreasing and that $u^*$ is the generalized right-continuous inverse of $\mu_u$.  Analogously one can define an increasing rearrangement $u_*$ of a function $u$.

An important property of decreasing rearrangements is that $u$ and $u^*$ have the same distribution function, and therefore
\begin{equation}\label{ext-e36}
    \int_B f(u(\xi))\, d\xi
    = \int_0^\infty f(u^*(s))\,ds\quad\text{for all measurable $f:\real\to\real_+$.}
\end{equation}

Let $X$ be a L\'evy process with characteristic exponent $\psi$. Set  $u(\xi):=e^{-t\Re\psi(\xi)}$ and denote by
\begin{equation}\label{ext-e38}
    \nu_{\Re\psi}(s) := m\{\xi\::\: \Re\psi(\xi) \leq s\}.
\end{equation}
Then we find for the decreasing rearrangement of the function  $u$ that
\begin{equation}\label{ext-e40}
  \begin{split}
    u^*(s)
    &= \inf\big\{\tau \geq 0\::\: m\{\xi:\,\, u(\xi)>\tau \}<s \big\}\\
    &=\inf\left\{\tau\geq 0\::\: \nu_{\Re \psi}\left(-t^{-1}\ln\tau\right)<s\right\}\\
    &= \exp\left(-t\nu_{\Re \psi}^{-1}(s)\right), \quad s>0.
  \end{split}
\end{equation}
Here $\nu_{\Re \psi}^{-1}$ is the generalized right-continuous inverse of $\nu_{\Re \psi}$, i.e.\ the increasing rearrangement $(\Re\psi)_*$ of $\Re\psi$. In particular, $(\Re\psi)_*$ is increasing. This allows to apply the same arguments as in the proof of Theorem~\ref{main-03} to find a necessary and sufficient criterion when $e^{-t\psi}\in L_1(\rn)$. It is worth noticing that the increasing rearrangement has an geometric meaning since it is the inverse volume growth function of metric balls of the form $\{\xi\in\rd \::\: d_\psi(\xi,0)\leq r\}$ with the intrinsic metric $d_\psi(\xi,\eta):=\sqrt{\psi(\xi-\eta)}$---see also Lemma \ref{ext-17} and Proposition \ref{app-21} below.
\begin{prop}\label{ext-11}
Let $X$ be an $n$-dimensional L\'evy process with characteristic exponent $\psi$. Then the following conditions are
equivalent
\begin{enumerate}
    \item[(a)] $e^{-t\psi}\in L_1(\rn)$ for all $t>0$;
    \item[(b)] $p_t(dx)=p_t(x)\,dx$ for all $t>0$, and $p_t\in  L_1(\rn)\cap C_\infty(\rn)$;
    \item[(c)] $p_t(dx)=p_t(x)\,dx$ for all $t>0$, and $p_t\in  L_\infty(\rn)$;
    \item[(d)] The increasing rearrangement $(\Re \psi)_*$ satisfies the following Hartman-Wint\-ner-type condition:
        \begin{equation}\label{hw-type}
            \lim_{|\xi|\to\infty} \frac{(\Re \psi)_*(\xi) }{\ln (1+|\xi|)} =\infty.
            \tag{$\textup{HW}^*_\infty$}
        \end{equation}
\end{enumerate}
\end{prop}
\begin{proof}
(a)$\Rightarrow$(b): By the Riemann-Lebesgue Lemma we see that $p_t = \Ff e^{-t\psi} \in C_\infty(\rn)$; since $p_t$ is the density of a probability measure, $p_t\in L_1(\rn)$ is automatically satisfied.

\medskip\noindent
(b)$\Rightarrow$(c): This is obvious.

\medskip\noindent
(c)$\Rightarrow$(a):
Since $p_t$ is a transition density, $p_{t/2}\in L_\infty (\rn) \cap L_1(\rn)\subset L_2(\rn)$, which implies
\begin{gather*}
    e^{-(t/2)\psi} = \Ff^{-1} p_{t/2} \in L_2(\rn)
    \implies
    |e^{-t\psi}| \in L_1(\rn)
    \implies e^{-t\psi}\in L_1(\rn).
\end{gather*}

\medskip\noindent
(d)$\Rightarrow$(a): This follows from \eqref{ext-e36} and \eqref{ext-e40}.

\medskip\noindent
(a)$\Rightarrow$(d): This follows from  \eqref{ext-e36}, \eqref{ext-e40} and the proof of (f)$\Rightarrow$(a) in Theorem~\ref{main-03}.
\end{proof}

The conditions \eqref{hw} and \eqref{hw1t} can be seen as comparison conditions: indeed we compare the growth rates (as $|\xi|\to\infty$) of the the characteristic exponent $\psi$ of the L\'evy process $X$ with the logarithm of the characteristic exponent of the symmetric Cauchy process. It is, therefore, a natural question how one can generalize  Theorem~\ref{main-03}.  Consider the following Hartman-Wintner-type condition
\begin{equation}\label{HW-new}\tag{$\textup{HW}^\phi_\infty$}
    \lim_{|\xi|\to\infty} \frac{\Re\psi(\xi)}{\ln(1+\phi(\xi))} = \infty
\end{equation}
where $\phi:\rn\to\real$ is the characteristic exponent of a further L\'evy process. For $\phi(\xi)=|\xi|$ and the corresponding Cauchy process \eqref{HW-new} and \eqref{hw} coincide. We are interested in the question which properties of $\phi$ imply that the L\'evy process $X$ with characteristic exponent $\psi$ admits a transition density $p_t$ such that $p_t\in C_\infty(\rn)\cap C^\infty(\rn)$.

For this we need to introduce a class of function spaces. Let $\phi$ be a real-valued characteristic exponent of a L\'evy process. Then $w = 1+\phi$ is a temperate weight function in the sense of \cite{hor} and it is possible to define the following class of \emph{anisotropic $L_2$-Bessel potential spaces}
\begin{equation*}
    H_2^{\phi,\kappa}(\rn)
    :=\left\{u\in S'(\mathbb{R}^n)\::\: \Ff^{-1} [(1+\phi)^{\kappa/2} \Ff u]\in L_2(\rn)\right\},\quad \kappa>0,
\end{equation*}
see \cite{FJS1,FJS2} for the general $L_p$-theory. For $p=2$ similar spaces were introduced by H\"{o}rmander \cite{hor} and, independently, by Volevich and Paneyah \cite{VP65} in order to study regularity properties of certain partial differential equations. For these spaces, the condition
\begin{equation}\label{ext-e42}
    \Ff^{-1}\left[(1+\phi)^{-\kappa/2}\right]\in L_2(\rn)
    \quad\text{or, equivalently,}\quad
    (1+\phi)^{-\kappa/2} \in L_2(\rn)
\end{equation}
is necessary and sufficient for the following Sobolev embedding
\begin{equation}\label{ext-e44}
    H_2^{\phi,\kappa}(\rn) \hookrightarrow C_\infty(\rn),
\end{equation}
see e.g.\ \cite[Theorem 2.3.4]{FJS2}.

In order to formulate the next theorem we need the notion of a Fourier integral resp.\ pseudodifferential operator: Let $\Psi:\rn\to\real$ be polynomially bounded. Then
$$
    \Psi(D) u(x) := \Ff^{-1}(\Psi\cdot \Ff u)(x) = \int_\rn e^{ix\xi}\,\Psi(\xi)\cdot\Ff u(\xi)\,d\xi,\quad u\in C_c^\infty(\rn),
$$
defines an linear operator with \emph{symbol} $\Psi$. If $\Psi(\xi)$ is a polynomial, $\Psi(D)$ is a differential operator, e.g.\ the symbol $\Psi(\xi)=\xi$ corresponds to the operator $D = -i\nabla$.

\begin{theorem}\label{ext-13}
    Let $X$ be an $n$-dimensional L\'evy process, $n\geq 1$, with characteristic exponent $\psi$ having no Gaussian part  and let $\phi$ be the characteristic exponent of a symmetric L\'evy process. If \eqref{ext-e42} holds for some $\kappa>0$, then the following conditions are equivalent:
    \begin{enumerate}
        \item[(a)] \eqref{HW-new};
        \item[(b)] for all $t>0$ the transition density exists, $p_t\in C_\infty(\rn)$ and $\phi(D)^m p_t\in L_1(\rn)\cap C_\infty(\rn)$ for all $m\in\nat$;
        \item[(c)] for all $t>0$ the transition density exists and $p_t, \phi(D)p_t \in L_1(\rn)$.
    \end{enumerate}
\end{theorem}
\begin{proof}
Denote by $\nu$ and $\mu$ the L\'evy measures of $\psi$ and $\phi$, respectively. As in the proof of Theorem \ref{main-03} we write $\psi_1$ and $\phi_1$ for the characteristic exponents which we obtain from $\Re\psi$ and $\phi$ by restricting the respective L\'evy measures to the unit ball: $\nu\rightsquigarrow \I_{\{0<|y|\leq 1\}}\,\nu(dy)$ and $\mu\rightsquigarrow \I_{\{0<|y|\leq 1\}}\,\mu(dy)$. Then we see that $\psi_1$ and $\phi_1$ are infinitely often differentiable and
$$
    \psi_1(\xi) \leq \Re \psi(\xi) \leq c(\psi_1(\xi) + 1)
    \quad\text{and}\quad
    \phi_1(\xi) \leq \phi(\xi) \leq c(\phi_1(\xi) + 1)
$$
for all $\xi\in\rn$ and some constant $c$. Thus, \eqref{HW-new} holds for $\psi$ if, and only if,  $(\textup{HW}^{\phi_1}_\infty)$ holds for $\psi_1$. Without loss of generality we may therefore assume that $\phi,\psi\in C^\infty(\rn)$.

\medskip\noindent
(a)$\Rightarrow$(b): Without loss of generality we may assume that $\psi$ is real-valued. By the chain rule
$$
    \nabla^\alpha e^{-t\psi} = (-t)^{|\alpha|}\,P(\psi,\ldots,\nabla^\alpha\psi) e^{-t\psi}
$$
where $P(\psi, \ldots, \nabla^\alpha\psi)$ is a polynomial of degree less or equal than $|\alpha|$ depending on the derivatives of $\psi$ of order $\beta\leq\alpha$. Using \eqref{proofs-e26}  and the Cauchy-Schwarz inequality  it is not hard to see that
$$
    \big|\nabla^\alpha \psi_1(\xi)\big|
    \leq
    c_{|\alpha|}
    \begin{cases}
    \displaystyle
        |\psi(\xi)|, &|\alpha|=0 \\[\medskipamount]
    \displaystyle
        \sqrt{\Re\psi(\xi)}, & |\alpha| = 1,\\[\medskipamount]
    \displaystyle
         1, & |\alpha| \geq 2.
    \end{cases}
$$
 For real-valued $\psi$  this observation is due to Hoh \cite{hoh-osaka}. The same is true for the exponent $\phi$. Therefore, $\xi\mapsto\nabla^\alpha e^{-t\psi(\xi)}$ is bounded for all $t>0$ and $\alpha\in\nat_0^n$. Combining this with \eqref{HW-new}, we find suitable polynomials $Q,R$ of degree less or equal than $|\alpha|+m$ such that for all $\kappa>0$ and $\alpha\in\nat_0^n$
\begin{align*}
    \left|(1+\phi)^\kappa \nabla_\xi^\alpha \big[(1+\phi)^m e^{-t\psi}\big](\xi)\right|
    &= \left|\sum_{\beta\leq\alpha} c_\beta \, \nabla^\beta (1+\phi)^m \nabla^{\alpha-\beta} e^{-t\psi}\right|\\
    &\leq c_t Q(\phi,\ldots,\nabla^\alpha\phi) R(\psi,\ldots,\nabla^\alpha\psi) e^{-t\Re\,\psi}
\end{align*}
is bounded. Because of \eqref{ext-e42},
$$
    \nabla_\xi^\alpha \big[(1+\phi)^m e^{-t\psi}\big](\xi)
    \in L_1(\rn)\cap L_\infty(\rn).
$$
From the Riemann-Lebesgue Lemma we get
$$
    |x|^k (1+\phi(D))^m p_t(x) \in C_\infty(\rn)
$$
for all $k,m\geq 0$. Choosing $k>n$ we can divide by $(1+|x|^k)$ and find, in particular, $(1+\phi(D))^m p_t(x)\in L_1(\rn)$.

\medskip\noindent
(b)$\Rightarrow$(c): This is obvious.

\medskip\noindent
(c)$\Rightarrow$(a): This follows as in the proof of Theorem \ref{main-03}.
\end{proof}

\begin{rem}\label{ext-15}
    \textbf{(a)} If we choose $\phi=\Re\psi$ in Theorem \ref{ext-13}, \eqref{HW-new} holds if, and only if, $\lim_{|\xi|\to\infty}\Re\psi(\xi)=\infty$.  In this case the statement
    \begin{enumerate}
    \item[(i)]
        \eqref{ext-e42} holds for $\phi=\Re\psi$ and some $\kappa>0$, i.e.\ $(1+\Re\psi)^{-\kappa/2}\in L_2(\rn)$
    \end{enumerate}
    always implies
    \begin{enumerate}
    \item[(ii)]
        for all $t>0$ we have $p_t(dx)=p_t(x)\,dx$ with $p_t\in \bigcap_{r>0} H^{\Re\psi,r}_2(\rn)$ and $\Re \psi(D)^m p_t\in C_\infty(\rn)$ for all $m\geq 0$.
    \end{enumerate}
    This follows follows immediately from Theorem \ref{ext-13}, \eqref{ext-e42} and \eqref{ext-e44}.

    We cannot expect the converse to be true. Even if  \eqref{HW-new} holds true, the condition $\phi^m(D)p_t\in  C_\infty(\rn)$ does not imply   \eqref{ext-e42}. To see this, consider the function $\psi$ from Example~\ref{exa2}. In this example $\psi\in C^\infty(R)$ and  $\psi(\xi)\sim \ln^2|\xi|$ as $|\xi|\to\infty$. Let $\phi(\xi):=\ln (1+\psi(\xi))$. Since $\psi$ and $\phi$ are infinitely often differentiable and, by $(HW_\infty)$, the function $\nabla^\alpha\big[\phi(\xi)e^{-t\psi(\xi)}\big]$ is  bounded for any $m,\alpha\geq 1$,  we have $\phi \cdot e^{-t\psi}\in S(R)$, implying
     $$
         \phi^m (D)p_t\in L_1(R)\cap C_\infty(R).
     $$
     But since $\phi(\xi)\sim \ln\ln |\xi|$ as $|\xi|\to\infty$, the condition   \eqref{ext-e42} does not hold for any $\kappa$.

    \bigskip\noindent
    \textbf{(b)} Theorem \ref{main-05} has also an obvious generalization: if we replace in the statement \eqref{hw} by \eqref{HW-new}, then we have to change $|\ln\varepsilon|$ to $\ln \phi(1/\varepsilon)$ in \eqref{main-e06}. The proof should be fairly obvious.
\end{rem}

The proof of Theorem \ref{ext-13} shows that the embedding condition \eqref{ext-e42} is important to assure that $p_t$ and $\phi(D)^m p_t$ are continuous functions. Let us  briefly discuss the meaning of \eqref{ext-e42}.
\begin{lem}\label{ext-17}
    Let $\phi$ be a real-valued characteristic exponent of some $n$-dimensional L\'evy process such that $\lim_{|\xi|\to\infty}\phi(\xi)=\infty$. Then \eqref{ext-e42} is equivalent to
    \begin{gather}\tag{$\ref{ext-e42}'$}\label{ext-e42*}
        \nu_\phi(x) := m(\xi\in\rn \::\: \phi(\xi) \leq x) \leq c\,x^\lambda
        \quad\text{for all $x\geq 1$}
    \end{gather}
    where $c,\lambda>0$ are suitable constants.

    The condition \eqref{ext-e42*}, in turn, is implied by the following volume-doubling condition
    \begin{equation}\label{voodoo}
        \frac{\nu_\phi(2x)}{\nu_\phi(x)} \leq C
        \quad\text{for some constant $C>0$ and all $x\geq 1$}.
    \end{equation}
\end{lem}
\begin{proof}
    Recall that $\phi_*(x) = \nu_\phi^{-1}(x)$ is the increasing rearrangement of the function $\phi$. Therefore,
    \begin{equation}\label{ext-e48}
        \int_\rn \frac {d\xi}{(1+\phi(\xi))^\kappa}
        = \int_0^\infty \frac {dx}{(1+\phi_*(x))^\kappa}
        = \kappa \int_0^\infty \frac{1}{(1+x)^{\kappa+1}}\,\nu_\phi(x)\,dx;
    \end{equation}
    the last equality follows from integration by parts.

    If \eqref{ext-e42} holds, the above integrals are finite. Since $\phi_*$ is increasing, the usual Abelian trick (cf.\ the proof of (f)$\Rightarrow$(a) in Theorem \ref{main-03}) guarantees that
    $$
        \frac{1}{(1+\phi_*(x))^\kappa}\leq \frac c{1+x}
        \implies
        \frac{1}{(1+y)^\kappa}\leq \frac c{1+\nu_\phi(y)}
    $$
    which gives \eqref{ext-e42*}.

    Conversely, if \eqref{ext-e42*} holds, the second equality in \eqref{ext-e48} shows that the integrals are finite as soon as $\kappa>\lambda$. This gives \eqref{ext-e42}.

    Assume finally that \eqref{voodoo} holds. Fix $x\geq 1$ and set $k:= [\ln x / \ln 2]$. Since $\nu_{\phi}$ is increasing, we find
    $$
        \nu_{\phi}(x)
        \leq \nu_{\phi}(2^{k+1})
        \leq C^{k+1}\nu_{\phi}(1)
        = C\,x^{\lambda}\nu_{\phi}(1)
    $$
    with $\lambda = \ln C/\ln 2$. Without loss of generality we can assume that $\lambda > 0$, and \eqref{voodoo} follows.
\end{proof}

\section{Applications}

The integral representation \eqref{main-e22} allows to determine in some cases the asymptotic behaviour of $p_t(0)$ in terms of the L\'evy exponent $\psi$. In what follows we do not assume isotropy.

Write $\nu_{\Re\psi}(x) := m\left\{\xi \::\: \Re\psi(\xi)\leq x\right\}$ for the distribution function \eqref{ext-e38} of $\Re\psi$. Then
\begin{equation*}
    (2\pi)^n\, p_t(0)
    =\|e^{-t\psi}\|_{L_1}
    =\|e^{-t\Re\psi}\|_{L_1}
    =t\int_0^\infty \nu_{\Re\psi}(x)\,e^{-tx}\, dx.
\end{equation*}
The following proposition below is essentially Theorem 4 from \cite[Chapter XIII.5]{Fel}.

\begin{prop}\label{app-19}
Suppose that \eqref{hw} holds true  and  $\nu_{\Re\psi}(x)\sim x^{\rho-1}L(x)$, $1<\rho<\infty$, as $x\to\infty$ \textup{[}resp.\ as $x\to 0$\textup{]} where $L(x)$ is a function which is slowly varying at infinity \textup{[}resp.\ at zero\textup{]}. Then
\begin{equation*}
    p_t(0)\sim \frac{\Gamma(\rho)}{t^{\rho-1}}L(t^{-1})
    \quad\text{as\ \ }t\to 0\qquad [\text{resp.\ as\ \ }t\to\infty].
\end{equation*}
\end{prop}

The above statement can be generalized to the case when $\nu_{\Re\psi}$ is of bounded increase.
\begin{prop}\label{app-21}
    Suppose that \eqref{hw} holds true  and the  function $\nu_{\Re\psi}(x)$  satisfies  for $\lambda\geq 1$ the following volume doubling property:
    \begin{equation}\label{app-voodoo}
        \frac{\nu_{\Re\psi}(2 x)}{\nu_{\Re\psi}(x)}\leq C
    \quad\text{as\ \ }x\to \infty \qquad [\text{resp.\ as\ \ }x\to 0]
    \end{equation}
    for some constant  $C<\infty$. Then
    \begin{equation}\label{app-e52}
        c_1 \nu_{\Re\psi}\left(t^{-1}\right)
        \leq p_t(0)
        \leq c_2  \nu_{\Re\psi}\left(t^{-1}\right)
        \quad\text{for\ \ }t\to 0\qquad [\text{resp.\ for\ \ }t\to\infty].
    \end{equation}
\end{prop}
\begin{proof}
    Fix $\lambda\geq 1$ and set $k:= [\ln\lambda / \ln 2]$. Since $\nu_{\Re\psi}$ is increasing, we find
    $$
        \nu_{\Re\psi}(\lambda x)
        \leq \nu_{\Re\psi}(2^{k+1} x)
        \leq C^{k+1}\nu_{\Re\psi}(x)
        = C\lambda^{\alpha}\nu_{\Re\psi}(x)
    $$
    with $\alpha = \ln C/\ln 2$. Without loss of generality we can assume that $\alpha\geq 0$. Therefore, \eqref{app-voodoo} is equivalent to saying that
    \begin{equation}\label{app-e54}
        \frac{\nu_{\Re\psi}(\lambda x)}{\nu_{\Re\psi}(x)}\leq C(1+o(1)) \lambda^\alpha
    \quad\text{as\ \ }x\to \infty \qquad [\text{resp.\ as\ \ }x\to 0].
    \end{equation}

    It is enough to consider the case where $x\to\infty$, the arguments for $x\to 0$ are similar. Since $(2\pi)^n\,p_t(0)=\int_0^\infty \nu_{\Re\psi}(\frac{y}{t}) e^{-y}dy$ we have by monotonicity
    $$
        (2\pi)^n\,p_t(0)
        \geq \int_1^\infty \nu_{\Re\psi}(yt^{-1}) e^{-y}\,dy
        \geq \nu_{\Re\psi}(t^{-1}) \int_1^\infty e^{-y}\,dy.
        = c_1\, \nu_{\Re\psi}(t^{-1}) 
    $$
    Because of \eqref{app-e54},
    \begin{align*}
    (2\pi)^n \,p_t(0)
    &=\left(\int_0^1+\int_1^\infty \right)\nu_{\Re\psi}(yt^{-1}) e^{-y}\,dy \\
    &\leq \nu_{\Re\psi}(t^{-1}) \int_0^1e^{-y}dy + C\, \nu_{\Re\psi}(t^{-1}) \int_1^\infty y^\alpha e^{-y}\,dy\\
    &= c_2 \, \nu_{\Re\psi}(t^{-1}) 
    \qedhere
    \end{align*}
\end{proof}

\begin{rem}
    The volume doubling property \eqref{app-voodoo} for $\nu_{\Re\psi}$ is important: for example, if  $\psi(\xi)\sim \ln^2|\xi|$ as $|\xi|\to\infty$, one can show by the Laplace method, see \cite{Co65}, that $c_1\, e^{c_2/t}\leq p_t(0)\leq c_3 \, e^{c_4/t}$ for $t\in (0,1]$. We refer to \cite{KK10} for more results on transition density estimates in small time.
\end{rem}

\begin{rem}
    Let $\psi$ be a characteristic exponent of a L\'evy process. It is known that $\psi$ induces via $\rho(x,y):=\sqrt{\Re\psi(x-y)}$ a metric on $\rn$, see \cite[Lemma~3.6.21]{J01}. Define by $B(x,r;\rho):=\left\{y\in\rn\::\: \sqrt{\Re\psi(x-y)}\leq r\right\}$ a ball of radius $r$ centred at $x$ in the metric $\rho$. Then $\nu_{\Re\psi}(r)= m(B(x,\sqrt{r};\rho))$. This allows the following interpretation of Proposition \ref{app-21}: if the measure of a ball in the metric $\rho$ is regular enough, its behaviour at infinity controls the behaviour of the transition density at zero.
\end{rem}

\begin{exa}\label{exa5}
    Consider the isotropic case, i.e.\  $\psi(\xi)=g(|\xi|^2)$ for some continuous $g$, which we assume in addition to be monotone. Under the conditions of Proposition~\ref{app-21} we get for some $c_1, c_2>0$
    \begin{equation*}
        c_1 \big(g^{-1}(1/t)\big)^{\frac{n}{2}}
        \leq p_t(0)
        \leq c_2 \big(g^{-1}(1/t)\big)^{\frac{n}{2}},
    \end{equation*}
    as $t\to 0$  (resp.,  $t\to \infty$). As an application of  Proposition \ref{app-21} we note that \eqref{app-e52} immediately tells us that the asymptotic behaviour of the transition density of the tempered $\alpha$-stable or truncated $\alpha$-stable process is
    $$
        c_1 t^{-\frac{n}{\alpha}} \leq p_t(0)\leq c_2   t^{-\frac{n}{\alpha}}
        \quad\text{as $t\to 0$}.
    $$
    Indeed, since the real part of the characteristic exponent of a tempered $\alpha$-stable processes behaves like $Re\, \psi(\xi)\sim  c |\xi|^\alpha$, see \cite[Theorem 2.9]{R07}), we can apply Proposition~\ref{app-19} to get the asymptotic behaviour of $p_t(0)$. For the truncated $\alpha$-stable process the asymptotic behaviour of $p_t(0)$ as $t\to 0$ follows from Proposition~\ref{app-21} and the observation that the characteristic exponent $\psi_R$ behaves like  $\psi_R(\xi)\sim |\xi|^\alpha$ if $\psi(\xi)\sim   |\xi|^\alpha$ as $|\xi|\to\infty$.
\end{exa}

Let us finally prove a straighforward ratio-limit theorem for L\'evy processes. We begin with an approximation result.
\begin{lem}\label{app-27}
    Let $\psi:\rn\to\comp$ be the characteristic exponent of a L\'evy process given by \eqref{intro-e02}. Assume that $e^{-t\psi}\in L^1(\rn)$ for all $t\geq t_0$. Then the normalized function $\chi_t := e^{-t\psi}/\| e^{-t\psi}\|_{L^1}$ satisfies for all $\delta>0$
    \begin{equation}\label{app-e56}
        \lim_{t\to\infty} \int_{|\xi| > \delta} |\chi_t(\xi)|\,d\xi = 0.
    \end{equation}
\end{lem}
\begin{proof}
    Let $t>t_0$. Then $p_t(dx)=p_t(x)\,dx$ with $p_t\in L^1(\rn)$ and, by the Riemann-Lebesgue Lemma, $\Ff^{-1} p_t = e^{-t\psi}\in C_\infty(\rn)$. In particular, $\lim_{|\xi|\to\infty} \Re\psi(\xi)=\infty$ which means that the following infimum
    $$
        m_\delta := \inf_{|\xi|>\delta} \Re\psi(\xi) > 0,\quad\delta>0,
    $$
    is attained and strictly positive. Otherwise there would be some $\xi_0\neq 0$ with $\psi(\xi_0)=0$ and $\psi$ would be periodic; this follows at once from an inequality for the characteristic exponent $\psi$ of a L\'evy process
    $$
        \big|\psi(\xi)+\overline{\psi(\eta)}-\psi(\xi-\eta)\big|^2
        \leq 4 |\psi(\xi)| \, |\psi(\eta)|
        \quad\text{for all}\quad \xi,\eta\in\rn,
    $$
    cf.\ Jacob \cite[Lemma 3.6.21]{J01} or Berg and Forst \cite[p.~46, proof of Lemma 7.15]{ber-for}. In  this case we would find for all $\epsilon>0$ some $h=h_\epsilon$ such that for all $k\in\nat$ $e^{-t\psi}\big|_{B_h(k\xi_0)} > 1-\epsilon$. This, however would contradict the assumption that $e^{-t\psi}\in L^1(\rn)$.

    Moreover, since $\psi$ is unbounded, $\nu(B_\delta(0))=\infty$ for any $\delta>0$.

    Therefore, for all $t>t_0$
    \begin{align*}
        \int_{|\xi|>\delta} |e^{-t\psi(\xi)}|\,d\xi
        &= \int_{|\xi|>\delta} e^{-t\Re\psi(\xi)}\,d\xi\\
        &= \int_{|\xi|>\delta} e^{-(t-t_0)\Re \psi (\xi)}e^{-t_0\Re\psi(\xi)}\,d\xi\\
        &\leq e^{-(t-t_0)m_\delta}\int_{|\xi|>\delta} e^{-t_0\Re\psi(\xi)}\,d\xi.
    \end{align*}
    From the L\'evy-Khinchine formula \eqref{intro-e02} we get for every $R>0$
    $$
        \Re\psi(\xi)
        \leq c^\psi_R |\xi|^2 + d^\psi_R,
        \quad\xi\in\rn
    $$
    where $c^\psi_R \asymp \|Q\| + \int_{|y|\leq R}|y|^2\,\nu(dy)$ and $d^\psi_R\asymp \nu(B_R^c(0))$. Thus,
    \begin{align*}
        \int_{|\xi|>\delta} |\chi_t(\xi)|\,d\xi
        &\leq \frac{e^{-tm_\delta}\, e^{t_0 m_\delta} \int_{|\xi|>\delta} e^{-t_0\Re\psi(\xi)}\,d\xi}{\int_\rn e^{-t\Re \psi (\xi)}\,d\xi}\\
        &\leq \frac{e^{-t m_\delta}\, e^{t_0 m_\delta} \int_{|\xi|>\delta} e^{-t_0\Re\psi(\xi)}\,d\xi}{\int_\rn e^{-t c^\psi_R |\xi|^2}\,d\xi\, e^{-t d^\psi_R}}\\
        &= (\sqrt t)^n\, e^{-t(m_\delta - d^\psi_R)}\,\frac{e^{t_0 m_\delta} \int_{|\xi|>\delta} e^{-t_0\Re\psi(\xi)}\,d\xi}{\int_\rn e^{- c^\psi_R |\xi|^2}\,d\xi}.
    \end{align*}
    Now we choose $R$ so large that $m_\delta > d^\psi_R$ and let $t\to\infty$. This proves \eqref{app-e56}.
\end{proof}

The following result is, in one dimension, due to W.\ Schenk \cite{schenk}.
\begin{theorem}\label{app-29}
    Let $X_t$ be a L\'evy process in $\rn$ with characteristic exponent $\psi$ and with transition semigroup $T_t$. If $e^{-t\psi}\in L^1(\rn)$ for all $t\geq t_0$, then the following limits exist locally uniformly for all $x\in\rn$, resp.\ $x,y\in\rn$
    \begin{align}
    \label{app-e60}
        \lim_{t\to\infty} \frac{T_t f(x)}{\|e^{-t\psi}\|_{L^1}}
        &= \frac 1{(2\pi)^n}\int_\rn f(z)\,dz\qquad\text{for all\ \ } f\in L^1(\rn).\\
        \lim_{t\to\infty} \frac{T_t f(x)}{T_{s+t}g(y)}
    \label{app-e62}
        &= \frac{\int_\rn f(z)\,dz}{\int_\rn g(z)\,dz}\qquad\text{for all\ \ } f,g\in L^1(\rn),\; s\geq 0;\\
    \label{app-e63}
        \lim_{t\to\infty} \frac{p_t(x)}{p_t(0)}
        &= 1.
    \end{align}
\end{theorem}
\begin{proof}
    It is clearly enough to prove \eqref{app-e60}. For $u\in C_\infty(\rn)$ we get from Lemma \ref{app-27} with $\chi_t(\xi) = e^{-t\psi(\xi)}/\|e^{-t\psi}\|_{L^1}$
    $$
        \lim_{t\to\infty} \int_\rn \chi_t(\xi)\,u(\xi)\,d\xi = u(0).
    $$
    Indeed, $\int_\rn \chi_t(\xi)\,d\xi = 1$ and
    \begin{align*}
        \left|\int_\rn \chi_t(\xi)(u(\xi)-u(0))\,d\xi\right|
        &\leq \int_{|\xi|\leq\delta} |\chi_t(\xi)| |u(\xi)-u(0)|\,d\xi + 2\int_{|\xi|>\delta} |\chi_t(\xi)|\,d\xi \|u\|_\infty\\
        &\leq \sup_{|\xi|\leq\delta} |u(\xi)-u(0)|\cdot \|\chi_t\|_{L^1} + 2\int_{|\xi|>\delta} |\chi_t(\xi)|\,d\xi \|u\|_\infty.
    \end{align*}
    Because of \eqref{app-e56} the second term vanishes as $t\to\infty$. Letting $\delta\to 0$ makes the first term tend to zero since $u$ is uniformly continuous.

    For $f\in L^1(\rn)$ and $t>t_0$ we have
    $$
        \frac{T_t f(x)}{\|e^{-t\psi}\|_{L^1}}
        = \frac 1{\|e^{-t\psi}\|_{L^1}}\,\Ff^{-1}\left[e^{-t\psi}\Ff f\right](x)
        = \int_\rn e^{ix\xi} \chi_t(\xi) \Ff f(\xi)\,d\xi.
    $$
    The above calculation shows for $u(\xi) := e^{ix\xi}\,\Ff f(\xi)$ and uniformly for $x$ from compact sets that
    \begin{gather*}
        \frac{T_t f(x)}{\|e^{-t\psi}\|_{L^1}}
        \xrightarrow{t\to\infty}
        \Ff f(0)
        = (2\pi)^{-n} \int_\rn f(z)\,dz.
    \qedhere
    \end{gather*}
\end{proof}

If we combine \eqref{app-e60} of Theorem \ref{app-29} with Propositions \ref{app-19}, \ref{app-21} or with Example \ref{exa5}, it is possible to get estimates for the speed of convergence in \eqref{app-e60}

\begin{ack}
    Financial support through the Ministry of Science of Ukraine (grant no.\ M/7-2008---for V.K.) and DFG (grant Schi 419/5-1---for R.L.S.) is gratefully acknowledged. We would like to thank an anonymous referee for his careful reading and helpful comments.
\end{ack}

\end{document}